\documentclass[a4paper]{amsart}
\usepackage[T1]{fontenc}
\usepackage{amssymb}
\usepackage[all]{xy}
\usepackage{verbatim}

\newcommand{\fun}[3]{#1\colon #2\rightarrow #3}

\newcommand{\set}[2]{\{#1:\, #2\}}
\newcommand{\bset}[2]{\Biggl\{#1:\, #2\Biggr\}}
\newcommand{\bbset}[2]{\Bigl\{#1:\, #2\Bigr\}}

\newcommand{\psr}[2]{#1[[#2]]}

\newcommand{\sepcl}[1]{\overline{#1}}
\newcommand{\algcl}[1]{\overline{#1}}

\newcommand{\abs}[1]{\lvert#1\rvert}
\newcommand{\absb}[1]{\bigl\lvert#1\bigr\rvert}

\newcommand{\cc}[1]{\mathcal{#1}}
\newcommand{\bb}[1]{\mathbb{#1}}

\newcommand{\dotsk}[0]{,\dots,}

\newcommand{\fG}[0]{\sepcl{\f}-\G}

\newcommand{\f}[0]{k}                     %Kropp
\newcommand{\sepclf}[0]{\sepcl{\f}}

\newcommand{\bbFq}[0]{\mathbb{F}_q}

\newcommand{\bbF}[0]{\mathbb{F}}
\newcommand{\bbL}[0]{\mathbb{L}}

\newcommand{\bbZ}[0]{\mathbb{Z}}
\newcommand{\bbA}[0]{\mathbb{A}}
\newcommand{\bbQ}[0]{\mathbb{Q}}
\newcommand{\multgroup}[0]{\mathbb{G}_m}          %Mult alg gruppen
\newcommand{\addgroup}[0]{\mathbb{G}_a}          %add alg gruppen

\newcommand{\ccP}[0]{\mathcal{P}}

\newcommand{\artcl}[0]{ArtCl}        %delring av artin klasser
\newcommand{\schur}[0]{Schur}

\newcommand{\catvar}[0]{\mathsf{Var}}                %K   "        variteer
\newcommand{\catgset}[1]{#1-\mathsf{Sets}}              %Kategorin av G-mängder
\newcommand{\catrep}[2]{#1-#2}

\newcommand{\per}[0]{\Sigma}                          %permutationsgrupp

\DeclareMathOperator{\KO}{K_0}

\DeclareMathOperator{\spec}{Spec}

\DeclareMathOperator{\Hom}{Hom}

\DeclareMathOperator{\V}{V}                             %Vershiebung
\DeclareMathOperator{\F}{F}                             %Frobenius

\DeclareMathOperator{\N}{N}

\DeclareMathOperator{\disjunion}{\overset{\centerdot}{\cup}}
\DeclareMathOperator{\Disjunion}{\overset{\centerdot}{\bigcup}}

\DeclareMathOperator{\D}{D}                             %schemaoperationer

\DeclareMathOperator{\Sym}{S}                             %symetric algebra
\DeclareMathOperator{\ext}{\bigwedge}                             %exterior algebra

\DeclareMathOperator{\Gal}{\mathcal{G}al}                 %Galoisgrupp
\DeclareMathOperator{\G}{\mathcal{G}}

\DeclareMathOperator{\image}{Im}

\DeclareMathOperator{\kar}{char}                          %karaktäristik

\DeclareMathOperator{\rep}{R}                          %representationsring
\newcommand{\ifc}[1]{C_{#1}}                           %karaktärhomomorfin
\newcommand{\karcyclo}[0]{\chi_{cycl}}                 %cyklotomisk karaktär

\DeclareMathOperator{\bur}{\mathcal{B}}                %Burnsidering

\DeclareMathOperator{\btr}{h}                        %burnside till representaionsring
\DeclareMathOperator{\res}{res}                        %restriction
\DeclareMathOperator{\ind}{ind}                        %induction

\DeclareMathOperator{\fro}{F}                         % Frobenius

\newcommand{\ress}[2]{R_{#1}^{#2}}                       %restriktion av skalärer
\newcommand{\inds}[2]{E_{#1}^{#2}}                     %utv av skaläre
\newcommand{\con}[1]{C_{#1}}                            %punkträkningshomomorfi
\DeclareMathOperator{\art}{Art}                        %artin avb.

\DeclareMathOperator{\aut}{Aut}
\DeclareMathOperator{\tr}{Tr}

\DeclareMathOperator{\homology}{H}

\numberwithin{equation}{section}

\newtheorem{theorem}[equation]{Theorem}
\newtheorem{lemma}[equation]{Lemma}

\newtheorem{proposition}[equation]{Proposition}
\newtheorem{example}[equation]{Example}
\newtheorem{definition}[equation]{Definition}
\newtheorem{remark}[equation]{Remark}

\theoremstyle{remark}
\newtheorem*{acknowledgements}{Acknowledgment}

\begin{document}

\title[The computation of the classes of some tori]{The computation of the classes of some tori in the Grothendieck ring of varieties}
\author{Karl Rökaeus}
\address{Karl Rökaeus \\ Matematiska institutionen \\ Stockholms universitet \\ SE-106 91 Stockholm \\ Sweden}
\email{karlr@math.su.se}
\date{August 31, 2007}
\begin{abstract}
In this paper we establish a formula for the classes of certain tori in the Grothendieck ring of varieties $\KO(\catvar_{\f})$. More explicitly, $\KO(\catvar_{\f})$ has a natural structure of a $\lambda$-ring, and we will see that if $L^*$ is the torus of invertible elements in the $n$-dimensional separable $\f$-algebra $L$ then $[L^*]=\sum_{i=0}^n(-1)^i\lambda^i\bigl([\spec L]\bigr)\mathbb{L}^{n-i}$, where $\mathbb{L}$ is the class of the affine line. This formula is suggested by the computation of the cohomology of the torus. To prove it will require some rather explicit calculations in $\KO(\catvar_{\f})$. To be able to make these we introduce a homomorphism from the Burnside ring of the absolute Galois group of $\f$, to $\KO(\catvar_{\f})$. In the process we obtain some information about the structure of the subring of $\KO(\catvar_{\f})$ generated by zero-dimensional varieties.
\end{abstract}

\maketitle

\section{Introduction}

The Grothendieck ring of varieties over the field $\f$, $\KO(\catvar_{\f})$, is the free abelian group on the objects of the category of varieties, subject to so called scissor relations and with a multiplication given by the product of varieties. Since the Euler characteristic with compact support (taking values for instance in a Grothendieck ring of mixed Hodge structures or Galois representations) is additive, it factors through $\KO(\catvar_{\f})$. Looking at relations among such Euler characteristics is a powerful heuristic method for guessing relations in $\KO(\catvar_{\f})$.

Our aim is to establish a formula for the classes in $\KO(\catvar_{\f})$ of certain tori, that can be guessed at in this manner. The cohomology of a torus can be expressed in terms of exterior powers of the first cohomology group. The latter in turn is essentially the cocharacter group of the torus. As $\KO(\catvar_{\f})$ is a $\lambda$-ring and as the Euler characteristic $\chi_c$ is a $\lambda$-homomorphism one can try to lift the cohomology formulas to $\KO(\catvar_{\f})$. There is a problem, however, in that in general one cannot find an element in $\KO(\catvar_{\f})$ that maps to the cocharacter representation under $\chi_c$. However, when the torus is the group $L^*$ of units in a separable $\f$-algebra $L$, then $[\spec L]$ maps to the cocharacter representation. We are thus led to conjecture the formula
\begin{equation}\label{3}
[L^*]=\sum_{i=0}^n(-1)^i\lambda^i\bigl([\spec L]\bigr)\bbL^{n-i}\in\KO(\catvar_{\f}),
\end{equation}
where $n$ is the dimension of $L$. (The details of this heuristic argument is given in the beginning of Section \ref{8}.) The main objective of this article is to prove \eqref{3}.

The structure of the paper is as follows: In Section \ref{2} we give a review of some definitions and results in order to fix notation.

In Section \ref{16} we construct a $\lambda$-homomorphism $\art_{\f}\colon\bur(\G)\to\KO(\catvar_{\f})$, where $\bur(\G)$ is the Burnside ring of the absolute Galois group of $\f$. The image of this map is contained in the subring of $\KO(\catvar_{\f})$ generated by zero-dimensional schemes, which we call the subring of Artin classes and denote $\artcl_{\f}$. This morphism immediately helps us answer some questions about the structure of $\artcl_{\f}$. For example, when $\f$ is perfect with cyclic absolute Galois group we will see $\artcl_{\f}$ is free on on the classes of finite separable field extensions of $\f$. Then in Section \ref{8} we use $\art_{\f}$ to simplify some of the computations by moving them to $\bur(\G)$.

In Section \ref{8} we show that $[L^*]=\sum_{i=0}^na_i\bbL^{n-i}\in\KO(\catvar_{\f})$, where $a_i\in\artcl_{\f}$. We also give an explicit formula for the $a_i$ in terms of elements in $\bur(\G)$. To derive this formula we embed $L^*$ in $\widetilde{L}$, the affine space associated to $L$, and use induction relative to the complement of $L^*$.

In the paper \cite{rokaeusBur} we have, by purely combinatorial methods, obtained a universal formula for the $\lambda$-operations on the Burnside ring which, together with the formula obtained in Section \ref{8}, gives a proof of \eqref{3}.

Finally, Section \ref{12} gives a proof of \eqref{3} that avoids using the universal formula for $\lambda^i$. Instead it uses one simpler result from \cite{rokaeusBur}, together with point counting over finite fields.

\begin{acknowledgements}
The author is grateful to Professor Torsten Ekedahl for valuable discussions and comments on the manuscript.

Thanks are also due to David Bourqui for comments on the previous version of this paper, in particular concerning Remark \ref{25}.
\end{acknowledgements}

\section{Review of some results}\label{2}
First a general remark. Many of the rings that we work with are constructed as free abelian groups on the objects of a category, subject to some relations. When defining a map from such a ring we often just give its action on an object in the category. We then have to show that it respects the relations. Also, when letting such a map act on the class of an object we often leave out the brackets, e.g., if $f\colon\KO(\catvar_{\f})\to R$ and $X$ is a scheme we write $f(X)$ for $f\bigl([X]\bigr)$.

\subsection*{$\lambda$-rings}
For an introduction, see for example the first part of \cite{MR0244387} or \cite{MR0364425}. A $\lambda$-ring is a commutative ring $R$ together with a homomorphism $\lambda_t$ from the additive group of $R$ to the multiplicative group of $\psr{R}{t}$, taking $x\in R$ to $\sum_{n\geq0}\lambda^n(x)t^n$, where $\lambda^0(x)=1$ and  $\lambda^1(x)=x$. A morphism of $\lambda$-rings $R\to R'$ is a ring homomorphism that commutes with the $\lambda^n$. Informally, this definition ensures that the $\lambda^n$ behave like exterior powers. The archetypal example is the representation ring of a finite group $G$. In this ring, $\lambda^n(V)=[\ext^n V]$, the $n$th exterior power of the vector space $V$ with componentwise $G$-action.

Let $\sigma_t$ be a $\lambda$-structure on $R$. The opposite structure of $\sigma_t$ is defined implicitly by the relation $\sigma_t(x)\lambda_{-t}(x)=1\in\psr{R}{t}$. On the representation ring the natural $\lambda$-structure can be obtained as the opposite to the one coming from the symmetric powers.

\subsection*{Representation rings}
We use $\rep_{\f}(G)$ to denote the ring of $\f$-representations of $G$, where $G$ is a profinite group. We require such a representation to be finite and continuous with respect to the profinite topology on $G$ and the \emph{discrete} topology on $\f$.

When $\f$ has some natural topology we can use the same construction but with respect to this topology instead. We call the ring thus obtained the Grothendieck ring of $\f$-representations of $G$, and denote it $\KO(\catrep{\f}{G})$. We have an injection $\rep_{\f}(G)\to\KO(\catrep{\f}{G})$, but this is not an isomorphism in general. For example, the cyclotomic representation is often not discrete.

As abelian groups, both these rings are free on isomorphism classes of irreducible representations. They are naturally $\lambda$-rings, the structure being given by exterior powers. A map $H\to G$ gives rise to an induction and a restriction map between the corresponding representation rings, which we denote $\ind_H^G$ and $\res_H^G$ respectively. Finally, for $g\in G$, we use $\ifc{g}$ to denote the character homomorphism from any representation ring of $G$, i.e., the map given by $V\mapsto\chi_V(g)$. Together they can be used to distinguish elements in the representation ring; $\prod_{g\in G}\ifc{g}$ is injective.

\subsection*{The Grothendieck ring of varieties}
For an introduction, see \cite{MR1886763} and the references given there. Let $\catvar_{\f}$ be the category of varieties over the field $\f$. Then $\KO(\catvar_{\f})$ is the free abelian group generated by symbols $[X]$ for $X\in\catvar_{\f}$, subject to the relations that $[X]=[Y]$ if $X\simeq Y$, $[X]=[X\setminus Y]+[Y]$ if $Y$ is a closed subscheme of $X$, and with a multiplication given by $[X]\cdot[Y]:=[X\times_{\f}Y]$. The second relation is usually referred to as the \emph{scissor relation}.
\footnote{
With respect to the definition of $\KO(\catvar_{\f})$, the important characterization of a variety is that it is of finite type over the base field; if not we end up with the zero ring. If we also include in the definition of a variety that it be reduced, we get a canonically isomorphic ring, for every scheme $X$ has a closed subscheme $X_{red}$ that is reduced and with empty complement, hence $[X]=[X_{red}]$. In the same way one can  add the conditions that a variety be separated and irreducible and still get an isomorphic ring. However, the condition that a variety be geometrically reduced probably gives a slightly smaller ring when $\f$ is non-perfect.}
By the class of the $\f$-scheme $X$ we mean its image $[X]\in\KO(\catvar_{\f})$. The class of the affine line is called the \emph{Lefschetz class} and denoted by $\bbL$. For a quick example of how the relations work, consider the multiplicative group $\multgroup$. It can be embedded in the affine line and its complement is then $\spec\f$. Hence $[\multgroup]=\bbL-1\in\KO(\catvar_{\f})$.

The Euler characteristic gives a map to the Grothendieck ring of $\bbQ_l$-representations of the absolute Galois group of $\f$, $$\chi_c\colon\KO(\catvar_{\f})\to\KO(\catrep{\bbQ_l}{\G}),$$ such that $\chi_c(X)=\sum(-1)^i[\homology^i_c(X_{\sepclf},\bbQ_l)]$. (Here $l$ is a prime different from the characteristic of $\f$.) We will use $\chi_c$ to study $\KO(\catvar_{\f})$, for example, its existence immediately shows that $\bbZ\subset\KO(\catvar_{\f})$.

Define a $\lambda$-structure $\KO(\catvar_{\f})$ as the opposite structure of $\{\sigma^n\}$, where $\sigma^n(X)=[X^n/\per_n]$. In \cite{MR2098401} it is  mentioned that this seems to be the natural $\lambda$-structure on $\KO(\catvar_{\f})$, since $\sigma^n$ behaves like a symmetric power map. Hence the Euler characteristic $\chi_c$ is a $\lambda$-homomorphism.

Let $K$ be a field extension of $\f$. We write $\ress{\f}{K}\colon\KO(\catvar_{K})\to\KO(\catvar_{\f})$ for the map defined by restricting scalars from $K$ to $\f$, i.e., $\ress{\f}{K}$ is defined by the map $\catvar_{K}\to\KO(\catvar_{\f})$ that take the $K$-scheme $X$ to the class of $X$, viewed as a $\f$-scheme. $\ress{\f}{K}$ is additive but not multiplicative. Also, define $\inds{\f}{K}\colon\KO(\catvar_{\f})\to\KO(\catvar_{K})$ by mapping the $\f$-scheme $X$ to the class of $X\times_{\f}\spec K$, viewed as a $K$-scheme. This is a ring homomorphism.

The fact that $\ress{\f}{K}$ fails to be multiplicative is because it does not preserve the multiplicative identity element, instead $\ress{\f}{K}(1)=[\spec K]\in\KO(\catvar_{\f})$. Rather than being multiplicative, $\ress{\f}{K}$ has a similar property: If $X$ is $\f$-scheme and $Z$ is a $K$-scheme then, from the universal property defining fibre products, $Z\times_K(\spec K\times_{\f}X)\simeq Z\times_{\f}X$ as $\f$-schemes. It follows that if we apply the restriction map to $[Z]\cdot[X_K]\in\KO(\catvar_{K})$ we get $[Z]\cdot[X]\in\KO(\catvar_{\f})$, i.e., for $x\in\KO(\catvar_{\f})$ and $z\in\KO(\catvar_{K})$ we have the projection formula
\begin{equation*}
\ress{\f}{K}\bigl(z\cdot\inds{\f}{K}(x)\bigr)=\ress{\f}{K}(z)\cdot x.
\end{equation*}
We will use the special case when $X=\bb{A}^n_{\f}$ and $Z=\spec L$, where $L$ is a finite-dimensional $K$-algebra.
\begin{equation}\label{13}
\ress{\f}{K}([\spec L]\cdot\bbL^n)=[\spec L]\cdot\bbL^n\in\KO(\catvar_{\f}).
\end{equation}
In particular, $\ress{\f}{K}(\bbL^n)=[\spec K]\cdot\bbL^n$.

\subsection*{Burnside rings}
For an introduction to Burnside rings, as well as proofs of the statements below, see \cite{MR0364425} chapter II, 4.

If $G$ is a profinite group then $\catgset{G}$ is the category with objects finite sets with continuous $G$-actions (with respect to the inverse limit topology on $G$) and morphisms $G$-equivariant maps of such sets. We will denote the set of morphisms between the $G$-sets $S$ and $T$ by $\Hom_G(S,T)$. The Burnside ring of $G$, $\bur(G)$, is constructed from this category as the free abelian group generated by the symbols $[S]$, for every continuous $G$-set $S$, subject to the relations that $[S\disjunion T]=[S]+[T]$, that $[S]=[T]$ if $S\simeq T$, and with a multiplication given by $[S]\cdot[T]:=[S\times T]$, where $G$ acts diagonally on $S\times T$.

Since every $G$-set can be written as a disjoint union of transitive $G$-sets we see that the transitive sets generate $\bur(G)$, and in fact it is free on the isomorphism classes of these. Moreover, every finite transitive $G$-set is isomorphic to $G/H$ where $H$ is a subgroup, and $G/H\simeq G/H'$ if and only if $H$ and $H'$ are conjugate subgroups. So every element of $\bur(G)$ can be written uniquely as $\sum_{H\in R}a_H[G/H]$, where $R$ is a system of representatives of the set of conjugacy classes of subgroups of $G$ and where $a_H\in\bbZ$ for every $H$.

We give a $\lambda$-structure $\lambda_t$ on $\bur(G)$, by first defining $\sigma_t\colon\bur(G)\to\psr{\bur(G)}{t}$ by the map that takes the $G$-set $S$ to the power series $\sum_{n\geq0}[S^n/\per_n] t^n\in\psr{\bur(G)}{t}$, where $\per_n$ acts on $S^n$ by permuting the entries. We then define $\lambda_t$ as the structure opposite to $\sigma_t$. It is non-special, but should still be considered to be the natural $\lambda$-structure on $\bur(G)$.

We have a map to the rational representation ring, $\btr\colon\bur(G)\to\rep_{\bbQ}(G)$, which is defined by associating to the $G$-set $S$ the class of the permutation representation $\bbQ[S]$. This is a homomorphism of $\lambda$-rings, which is one of the reasons why we consider our $\lambda$-structure on $\bur(G)$ to be the right one. It has the property that if $S$ and $S'$ are two non-isomorphic transitive $G$-sets then $\btr(S)\neq\btr(S')$. However, since $\bur(G)$ has rank equal to the number of conjugacy classes of subgroups of $G$ whereas $\rep_{\bbQ}(G)$ has rank equal to the number of conjugacy classes of cyclic subgroups of $G$, $\btr$ cannot be injective unless $G$ is cyclic. Conversely, if $G$ is cyclic then $\btr$ is an isomorphism. (In this paper, when we say that a profinite group is cyclic we mean topologically, i.e, it has a dense cyclic subgroup.) These facts are proved for example by using the character maps $\ifc{g}$ for $g\in G$. (Usually for finite groups, but since taking the inverse limit of groups corresponds to taking the direct limit of the corresponding Burnside and representation rings, they follow immediately for any profinite group.)

Next, let $\phi\colon H\to G$ be a group homomorphism. We associate to it two maps, restriction and induction, between the corresponding Burnside rings in the same way as for representation rings: Firstly, $\res_H^G\colon\bur(G)\rightarrow\bur(H)$ is the map induced by restricting the $G$-action on a $G$-set $S$ to a $H$-action, i.e., $S$ is considered as a $H$-set via $h\cdot s:=\phi(h)s$ for $h\in H$ and $s\in S$. This map is a ring homomorphism.

Secondly, if instead $S$ is a $H$-set then we can associate to it the $G$-set $G\times_HS$, i.e., the quotient of $G\times S$ by the equivalence relation $(g\cdot\phi(h),s)\sim(g,hs)$ for $(g,s)\in G\times S$ and $h\in H$, with a $G$-action given by $g'\cdot(g,s):=(g'g,s)$. This gives rise to the induction map $\ind_H^G\colon\bur(H)\rightarrow\bur(G)$, which is additive but not multiplicative. We will only use the induction map in the case when $H$ is a subgroup of $G$. In this case, note that if we choose a set of coset representatives of $G/H$, $R=\{g_1\dotsk g_r\}$, then we can represent $G\times_HS$ as $R\times S$ with $G$-action given by $g\cdot(g_i,s)=(g_j,hs)$, where $gg_i=g_jh$ for $h\in H$.

The map $\btr\colon\bur(G)\to\rep_{\bbQ}(G)$ commutes with the restriction maps, and also with the induction maps if $H$ is a subgroup of $G$.

We conclude this subsection by giving an example of a Galois group $\G$ of a finite extension of $\bbQ$ such that $\btr_{\G}\colon\bur(\G)\to\rep_{\bbQ}(\G)$ is not surjective. (This will be used in Remark \ref{25}, to find tori for which our main theorem do not hold.) For this we use an example of Serre of a finite group having this property, together with the following lemma:
\begin{lemma}
Let $N$ be a normal subgroup of the finite group $G$, and let $H:=G/N$. If $x\in\rep_{\bbQ}(H)$ is not contained in the image of $\btr_H\colon\bur(H)\to\rep_{\bbQ}(H)$, then $\res_G^H x$ is not contained in the image of $\btr_G\colon\bur(G)\to\rep_{\bbQ}(G)$. In particular, if $\btr_G$ is surjective, then so is $\btr_H$.
\end{lemma}
\begin{proof}
Define a map $\rep_{\bbQ}(G)\to\rep_{\bbQ}(H)$ by $V\mapsto[V^N]$, where $V^N$ is the elements of the $G$-representation $V$ invariant under $N$. The corresponding map $\bur(G)\to\bur(H)$ is defined by $S\mapsto [S/N]$, and one proves that the following diagram is commutative:
$$
\xymatrix{
\bur(H)\ar[r]^{\res_G^H}\ar[d]^{\btr} & \bur(G)\ar[r]\ar[d]^{\btr} & \bur(H)\ar[d]^{\btr}\\
\rep_{\bbQ}(H)\ar[r]^{\res_G^H} & \rep_{\bbQ}(G)\ar[r] & \rep_{\bbQ}(H)
}
$$
Since both the horizontal compositions equal the identity, the result follows by diagram chasing.
\end{proof}
\begin{example}\label{24}
There exists a finite extension $\f$ of $\bbQ$, with absolute Galois group $\G$, such that $\btr_{\G}\colon\bur(\G)\to\rep_{\bbQ}(\G)$ is not surjective. For by Exercise 13.4, page 105 of \cite{MR0450380}, there is a finite group $H$ such that $\btr_H\colon\bur(H)\to\rep_{\bbQ}(H)$ is not surjective. (More precisely, the example in \cite{MR0450380} shows that this holds for the product of the quaternion group with the cyclic group of order $3$.) Now choose a finite field extension $\f$ of $\bbQ$ such that there exists a finite extension $K/\f$ with the property that $\Gal(K/\f)=H$. Choose $x\in\rep_{\bbQ}(H)$ that is not contained in the image of $\btr_H$. If $\btr_{\G}$ is surjective, then there is a finite quotient $G$ of $\G$, having $H$ as a quotient, such that $\res_G^Hx$ is in the image of $\btr_G$. By the preceding lemma, this is impossible.
\end{example}

\subsection*{Galois theory}
To be able to make explicit computations in $\KO(\catvar_{\f})$ we use Grothendieck's formulation of Galois theory. We use the following notation. Let $\f$ be a field and let $\sepclf$ be a separable closure of $\f$. Let $\G:=\Gal(\sepclf/\f)$, the absolute Galois group. Then the category of separable $\fG$-algebras is defined to be the category whose objects are separable $\sepclf$-algebras $L$ together with $\G$-actions on the underlying rings such that $\sepclf\rightarrow L$ is $\G$-equivariant, and whose morphisms are $\G$-equivariant maps of $\sepclf$-algebras.

We then have an equivalence between the category of $\f$-algebras and the category of $\fG$-algebras. This equivalence takes the $\f$-algebra $L$ to $L\otimes_\f\sepclf$ with $\G$-action $\sigma(l\otimes\alpha):=l\otimes\sigma(\alpha)$. Its pseudo-inverse takes the $\fG$-algebra $U$ to $U^{\G}$.

If the $\f$-algebra $L$ is finite and separable, then the corresponding $\fG$ algebra is isomorphic to $\sepclf^S$ where $S$ is finite. The $G$-action on this must be the action on $\sepclf$ together with a permutation of the coordinates, i.e., a $\G$-action on $S$. Hence as a corollary we have a contravariant equivalence between the category of finite separable $\f$-algebras and the category of finite continuous $\G$-sets. This equivalence takes the $\f$-algebra $L$ to $\Hom_{\f}(L,\sepclf)$ with $\G$-action given by $f^\sigma(l):=\sigma\circ f(l)$. Its pseudo-inverse takes the $\G$-set $S$ to $\Hom_{\G}(S,\sepclf)$, i.e., the $\G$-equivariant maps of sets from $S$ to $\sepclf$, considered as a ring by pointwise addition and multiplication and with a $\f$-algebra structure given by $(\alpha\cdot f)(s):=\alpha\cdot f(s)$.

Under this correspondence, if $L$ corresponds to $S$ then the dimension of $L$ equals the number of elements in $S$. Moreover, if also $L'$ corresponds to $S'$, then $L\otimes_{\f} L'$ corresponds to $S\times S'$ with diagonal $\G$-action and the algebra $L\times L'$ corresponds to $S\disjunion S'$. In particular, separable field extensions of $\f$ correspond to transitive $\G$-sets.

\section{The subring of Artin classes in $\KO(\catvar_{\f})$}\label{16}

In this section we define a map from the Burnside ring of the absolute Galois group of $\f$ to $\KO(\catvar_{\f})$. This map will be used in the computation of the class of $L^*$, but it also gives some information about the structure of $\KO(\catvar_{\f})$.
\begin{definition}
Let $\artcl_{\f}$, the ring of Artin classes, be the subring of $\KO(\catvar_{\f})$ generated by zero-dimensional schemes.
\end{definition}
As an abelian group $\artcl_{\f}$ is generated by $\{[\spec K]\}$, where $K$ runs over the isomorphism classes of finite field extensions of $\f$. When the characteristic of $\f$ is zero there is a structure theory for $\KO(\catvar_{\f})$ given in \cite{MR1996804}, which asserts that the classes of stable birational equivalence form a $\bbZ$-basis of $\KO(\catvar_{\f})/(\bbL)$. Using this, one shows that $\artcl_{\f}$ is free on $\{[\spec K]\}$. That $\artcl_{\f}$ is free on this set is also true when $\f=\bbFq$ as is shown in \cite{MR2272145}, Theorem 25, using the point counting homomorphisms $\con{q^n}$, where, for any prime-power $q$, $\con{q}\colon\KO(\catvar_{\bbF_q})\to\bbZ$ is given by $X\mapsto\abs{X(\bbF_q)}$.

There is also a question about zero divisors. In \cite{MR1928868} it is shown that when $\kar\f=0$, $\KO(\catvar_{\f})$ contains zero divisors. To show that these really are nonzero requires the structure theory of \cite{MR1996804}. \cite{MR2272145} observes that if $K$ is a finite Galois extension of degree $n$ then $[\spec K]^2=n\cdot[\spec K]$ so $[\spec K]\in\artcl_{\f}$ is a zero divisor if it is not equal to $0$ or $n$, and that this is the case when $\f=\bbFq$.

Let $\f$ be a field with absolute Galois group $\G$. In this section we use Galois theory to define a $\lambda$-homomorphism from  $\bur(\G)$ to $\KO(\catvar_{\f})$, whose image is contained in $\artcl_{\f}$. The main purpose of doing this is that it aids the computations in $\KO(\catvar_{\f})$. Also it allows us to give a slightly generalization of the above mentioned results of \cite{MR2272145}.

\begin{definition}
Let $\art_{\f}\colon\bur(\G)\to\KO(\catvar_{\f})$ be the $\lambda$-homomorphism defined by associating, to the $\G$-set $S$, the class of its image under the fully faithful, covariant Galois functor to $\catvar_{\f}$.
\end{definition}
Since open disjoint union is a special case of the relations in $\KO(\catvar_{\f})$ this is well defined. Moreover it is multiplicative since the multiplications in both rings come from the product in the respective categories. Finally, to see that it commutes with the $\lambda$-structures, note that it suffices to check this on the opposite structures, and $\sigma^n$ is defined by the same universal property in both rings (namely the $n$th symmetric power) so this is clear.

Since $\bur(\G)$ is free on isomorphism classes of transitive $\G$-sets, and $\art_{\f}$ maps these to isomorphism classes of separable field extensions of $\f$, it follows that $\image\art_{\f}$ is free on $\{[\spec K]\}$ if and only if $\art_{\f}$ is injective. If $\f$ is perfect, $\image\art_{\f}=\artcl_{\f}$, hence the same holds for $\artcl_{\f}$ in this case.

Recall that we use $\btr$ to denote the natural $\lambda$-homomorphism $\bur(\G)\to\rep_{\bbQ}(\G)$. We use $i$ for the injection $\rep_{\bbQ}(\G)\to\KO(\catrep{\bbQ_l}{\G})$. If $S$ is a $\G$-set and $L$ the corresponding separable $\f$-algebra then the cohomology of $\spec L$ is $\bbQ_l[S]$, hence we have the following commutative diagram of $\lambda$-rings:
\begin{equation}\label{17}
\xymatrix{
\bur(\G)\ar[r]^{\art_{\f}}\ar[d]^{\btr} & \KO(\catvar_{\f})\ar[d]^{\chi_c}\\
\rep_{\bbQ}(\G)\ar[r]^i & \KO(\catrep{\bbQ_l}{\G})
}
\end{equation}
In particular, when $\f=\bbFq$ and $\fro\in\G$ is the Frobenius automorphism, if $S$ is a $\G$-set corresponding to the variety $X$ then $\ifc{\fro}\circ\chi_c(X)=\abs{S^{\fro}}$ and also $\abs{X(\bbFq)}=\abs{\Hom_{\bbFq}(\spec\bbFq,X)}=\abs{\Hom_{\G}(\{\bullet\},S)}=\abs{S^{\fro}}$,
consequently
\begin{align}\label{20}
\con{q}=\ifc{\fro}\circ\chi_c && \text{on }\image\art_{\bbFq},
\end{align}
showing that the character maps generalize point counting, a fact that we will use in Section \ref{12}.

As a first application of the commutativity of \eqref{17} we note that if $L$ and $L'$ are two non-isomorphic separable field extensions of $\f$, i.e., they correspond to two non-isomorphic transitive $\G$-sets $S$ and $S'$, then $[\spec L]\neq[\spec L']$ and $[\spec L]\notin\bbZ$. For it suffices to show that this is not the case for their images under $\chi_c$, i.e., that $\btr(S)\neq\btr(S')$ and that $\btr(S)\notin\bbZ$, and this is a known property of $\btr$. In particular:
\begin{proposition}
For any field $\f$, the isomorphism classes of finite Galois extensions determine distinct zero divisors in $\KO(\catvar_{\f})$.
\end{proposition}
This type of argument cannot be used to prove that $\image\art_{\f}$ is free on $\{[\spec K]\}$, since in general there are non-isomorphic $\G$-sets $S$ and $S'$ such that $\btr(S)=\btr(S')$. The exception is when $\G$ is cyclic, for then $\btr$ is an isomorphism. We state this as a proposition.
\begin{proposition}
If the absolute Galois group of $\f$ is cyclic then $\art_{\f}$ is injective.
\end{proposition}
As a corollary we get a result from \cite{MR2272145}, that $\artcl_{\bbFq}$ is free on $\{[\spec K]\}$.

To summarize: We have a good understanding of the structure of $\artcl_{\f}$ when $\kar\f=0$, and also when $\f$ is perfect with cyclic Galois group.

\subsection*{The behavior of $\art_{\f}$ with respect to restriction of scalars}
We next study how $\art_{\f}$ behaves with respect to restriction of scalars. The following proposition is due to Grothendieck but we have not been able to find a reference so we include a proof for completeness.
\begin{proposition}\label{Pp:60}
Fix a field $\f$ together with a separable closure $\sepclf$ and let $\G:=\Gal(\sepclf/\f)$. Let $K$ be a finite field extension of $\f$ such that $K\subset\sepclf$. Let $L$ be a finite separable $K$-algebra and let $S$ be the corresponding $\Gal(\sepclf/K)$-set. View $L$ as a $\f$-algebra and let $S'$ be the corresponding $\G$-set. Then $S'\simeq\G\times_{\Gal(\sepclf/K)}S$. Hence the following diagram is commutative.
$$
\xymatrix{
\bur\bigl(\Gal(\sepclf/K)\bigr) \ar[r]^{\art_K}\ar[d]_{\ind_{\Gal(\sepclf/K)}^{\G}}  &  \KO(\catvar_{K}) \ar[d]^{\ress{\f}{K}}\\
\bur(\G)  \ar[r]^{\art_{\f}}  &   \KO(\catvar_{\f})
}
$$
\end{proposition}
\begin{proof}
Define a map $\phi\colon\G\times S\to S'$ by $(\sigma,f)\mapsto\sigma f$. It has the property that if $\tau\in\Gal(\sepclf/K)$ then $\phi(\sigma\tau, f)=\sigma\tau f=\phi(\sigma,\tau f)$. Hence it gives rise to a map of $\G$-sets $\varphi\colon\G\times_{\Gal(\sepclf/K)}S\to S'$. If $\phi(\sigma,f)=\phi(\tau,g)$ then $\tau^{-1}\sigma f=g$ so since $f$ and $g$ fixes $K$ pointwise we must have that $\tau^{-1}\sigma\in\Gal(\sepclf/K)$. It follows that $(\tau,g)=(\tau,\tau^{-1}\sigma f)\sim(\tau\tau^{-1}\sigma,f)=(\sigma,f)$ so $\varphi$ is injective. It is also surjective, for let $d:=[K:\f]$ and suppose that $L$ has dimension $n$ as a $K$-algebra, i.e., $S$ has $n$ elements. Then $L$ has dimension $nd$ as a $\f$-algebra so $S'$ has $nd$ elements. On the other hand, by Galois theory, $\abs{\G/\Gal(\sepclf/K)}=[K:\f]=d$. Hence $\G\times_{\Gal(\sepclf/K)}S$ also has $nd$ elements. Since $\varphi$ is injective it follows that it also is surjective, hence an isomorphism of $\G$-sets.
\end{proof}

\section{The class of a torus in $\KO(\catvar_{\f})$}\label{8}

Given a field $\f$ and a separable $\f$-algebra $L$ of dimension $n$ we define the affine group scheme $L^*$ by letting $L^*(R)=(L\otimes_{\f}R)^\times$ for every $\f$-algebra $R$. This is a torus, because if $\sepclf$ is a separable closure of $\f$ and $R$ is a $\sepclf$-algebra then, since $L\otimes_{\f}\sepclf=\sepclf^n$, we have $$L^*_{\sepcl{\f}}(R)=L^*(R)=\bigl((L\otimes_{\f}\sepcl{\f})\otimes_{\sepcl{\f}}R\bigr)^\times=(R^n)^\times=\multgroup^n(R)$$
as $\sepclf$-algebras, and consequently $(L^*)_{\sepclf}\simeq\multgroup^n$. We call $L^*$ the torus of invertible elements in $L$. Note that if $L=\f^n$ then $L^*=\multgroup^n$, hence $[L^*]=(\bbL-1)^n\in\KO(\catvar_{\f})$.

As mentioned in the introduction, if $T$ is a torus of dimension $n$ and $N$ is its cocharacter group tensored with $\bbQ_l$ then
\begin{equation}\label{19}
\chi_c(T)=\sum_{i=0}^n(-1)^i\lambda^i(N)\ell^{n-i}\in\KO(\catrep{\bbQ_l}{\G}),
\end{equation}
where $\ell:=[\bbQ_l(-1)]$, the class of the cyclotomic representation. If $T=L^*$, where $L$ corresponds to the $\G$-set $S$, then $[N]=[\bbQ_l[S]]=\chi_c(\spec L)$. Together with the fact that $\chi_c(\bbL)=\ell$, this suggests formula \eqref{3}.

The objective of this chapter is to compute, for an arbitrary separable $\f$-algebra $L$, the class of $L^*$ in $\KO(\catvar_{\f})$ in terms of the Lefschetz class $\bbL$ and Artin classes. More precisely, we will show that there exist elements $a_1\dotsk a_n\in\artcl_{\f}\subset\KO(\catvar_{\f})$ such that
\begin{equation}\label{21}
[L^*]=\bbL^n+a_1\bbL^{n-1}+a_2\bbL^{n-2}+\dots+a_n\in\KO(\catvar_{\f}).
\end{equation}
In Theorem \ref{111} we then derive a universal formula for the $a_i$ which, together with the result from \cite{rokaeusBur}, proves \eqref{3}.

Before we begin, let us mention that \eqref{21} does not hold for an arbitrary torus:
\begin{remark}\label{25}
There are fields $\f$ and tori $T$ such that $[T]\notin\artcl_{\f}[\bbL]$. To show this we need the following lemma: Identify $\rep_{\bbQ_l}(\G)$ with its image in $\KO(\catrep{\bbQ_l}{\G})$ and assume that $\abs{\image\karcyclo}=\infty$. If $\sum b_i\ell^i=0$, where $b_i\in\rep_{\bbQ_l}(\G)$, then $b_i=0$. This is proved using the character maps $\ifc{g}$, $g\in\G$, for we have that $\absb{\{\ifc{g}(\ell)\}_{g\in\G}}=\infty$ whereas $\absb{\{\ifc{g}(b_i)\}_{g\in\G}}<\infty$ for every $i$. Suppose now that $\f$ is such that $\abs{\image\karcyclo}=\infty$ and that $[T]\in\artcl_{\f}[\bbL]$, where $T$ is a torus of dimension $n$. Then, since $\chi_c(\bbL)=\ell$, we have $\chi_c(T)=\sum_{i=-m}^nb_i\ell^{n-i}$, where $m$ is an integer and $b_i$ is in $B$, the image of $\bur(\G)$ in $\KO(\catrep{\bbQ_l}{\G})$. Hence, by \eqref{19},
$\sum_{i=0}^n(-1)^i\lambda^i(N)\ell^{n-i}=\sum_{i=-m}^nb_i\ell^{n-i}.$
Consequently, since $\lambda^i(N)\in\rep_{\bbQ}(\G)\subset\rep_{\bbQ_l}(\G)$ and $b_i\in B\subset\rep_{\bbQ_l}(\G)$, the lemma shows that $b_i=(-1)^i\lambda^i(N)$. In particular, $[N]=\lambda^1(N)=-b_i\in B$. But not all tori has this property. To see that, note that every $\bbQ$-representation of $\G$, $V$, comes from the cocharacter group of a torus, namely the free $\bbZ$-module $\langle ge_i\rangle_{g\in\G}\subset V$, where $\{e_i\}$ is a basis for $V$. (Here it is important that the representation is discrete, otherwise this group can have infinite rank.) So every representation $[N]\in\rep_{\bbQ}(\G)$ corresponds to a torus $T$, and if $[N]$ is not contained in $B$, then $[T]\not\in\artcl_{\f}[\bbL]$. By Example \ref{24}, there is a finite extension of $\bbQ$, $\f$, such that the inclusion $B\subset\rep_{\bbQ}(\G)$ is proper. Since $\f$ is finitely generated, the image of the cyclotomic character is infinite, hence  it follows that there are tori $T$ such that $[T]\notin\artcl_{\f}[\bbL]$.
\end{remark}
\begin{remark}
The following was pointed out to us by David Bourqui: If $L$ is a separable $\f$-algebra of dimension $n$, and $S$ the corresponding $\G$-set, then the exact sequence of $G$-modules $0 \to \bbZ \to \bbZ[S] \to M\to 0$ splits over $\bbQ$, hence $[M\otimes_{\bbZ}\bbQ_l]=[\bbQ_l[S]]-1=\chi_c([\spec L]-1)$. So if $L^{*,1}$ is the torus with cocharacter group equal to $M$, i.e., $L^{*,1}=L^*/\multgroup$, then the above heuristic suggests that
\begin{equation}\label{26}
[L^{*,1}]=\sum_{i=0}^{n-1}(-1)^i\lambda^i\bigl([\spec L]-1\bigr)\bbL^{n-1-i},
\end{equation}
which also would follow from \eqref{3} if we knew that $\bbL-1$ was not a zero-divisor, because $[L^{*,1}]\cdot(\bbL-1)=[L^*]$. That \eqref{26} actually does hold can be proved by embedding $L^{*,1}$ in the projective space associated to $L $, and then use the exact same method as we will use to prove \eqref{3}.

We also mention one case when the above heuristic probably gives the wrong formula, namely let $T$ be the torus of elements in $L$ of norm $1$. If $M$ is the cocharacter group of $T$, then the exact sequence $0 \to M \to \bbZ[S] \to \bbZ\to 0$ splits over $\bbQ$, hence by the above argument, $[T]$ would satisfy the same formula as $[L^{*,1}]$. Looking at the top dimensional term of the formula suggests that $T$ should be rational, however, by \cite{MR0067122} page $322$, there are field extensions $L/\f$ such that $T$ is not rational.
\end{remark}

\subsection*{Definitions}
We begin by giving a definition of $L^*$ for any free algebra of finite rank: Let $K$ be a commutative ring and let $L$ be a free $K$-algebra of rank $n$. Let $\widetilde{L/K}$, or just $\widetilde{L}$, be $L$ considered as an affine space, i.e., the $R$-points on $\widetilde{L}$ is  $L\otimes_{K}R$. We have $\widetilde{L}=\spec\Sym(L^\vee)$. Note in particular that $\widetilde{K}$ is the ring scheme with additive group $\addgroup$ and multiplicative group $\multgroup$. Also, if we choose a $K$-basis for $L$ we get an isomorphism $\Sym(L^\vee)\simeq K[X_1\dotsk X_n]$, where $n$ is the rank of $L$. Hence $\widetilde{L}\simeq\bb{A}^n_{K}$ as schemes.

We next give the general definition of the $K$-scheme $(L/K)^*$, which we write as $L^*$ when $K$ is clear from the context. Define $L^*\subset\widetilde{L}$ as the subfunctor given by $L^*(R)=(L\otimes_{K} R)^\times$. To see that this is an affine group scheme, note that it is the inverse image of $K^*$ under the norm map $\fun{\widetilde{\N}_{L/K}}{\widetilde{L}}{\widetilde{K}}$. (Here, $\widetilde{\N}_{L/K}$ is defined on $R$-points as $\N_{L\otimes_{K}R/R}$.)

We now turn to the problem of computing $[(L/\f)^*]$ when $\f$ is a field and $L$ is separable. The obvious approach would be to compute an explicit equation defining $L^*$ and the use the scissor relations. More precisely, choose a basis of $L$ over $\f$. This identifies $\widetilde{L}$ with $\spec\f[X_1\dotsk X_n]$ and $\widetilde{\f}$ with $\spec\f[X]$. Then $\widetilde{\N}_{L/\f}$ corresponds to a homomorphism of $\f$-algebras $\f[X]\to\f[X_1\dotsk X_n]$ sending $X$ to the polynomial $f(X_1\dotsk X_n)$. Therefore, $\widetilde{\N}_{L/K}^{-1}(\bb{G}_m)=\spec\f[X,1/X]\otimes_{\f[X]}\f[X_1\dotsk X_n]=\spec\f[X_1\dotsk X_n,1/f(X_1\dotsk X_n)]$. When $n=2$ this can be used to compute $[L^*]$ using the scissor relations  as the following example shows.
\begin{example}\label{9}
Let $L$ be a separable extension field of $\f$ of degree $2$. We can represent $L$ as $\f[T]/\bigl(f(T)\bigr)$ where $f(T)=T^2+\alpha T+\beta$ is irreducible, in particular $\beta\neq0$. If $\kar\f\neq2$ we assume that $\alpha=0$. With this notation we have $\widetilde{L}(R)=R[T]/\bigl(f(T)\bigr)$ for every $\f$-algebra $R$. A basis for the $R$-algebra $\widetilde{L}(R)$ is $\{1,t\}$ where $t$ is the class of $T$ modulo $f(T)$. If $r_1,r_2\in R$ then $\N_{\widetilde{L}(R)/R}(r_1+r_2t)=r_1^2-r_1r_2\alpha+r_2^2\beta$.
So if we identify $\widetilde{L}$ with $\spec\f[X_1,X_2]$ then
$$L^*=\D\bigl(X_1^2-\alpha X_1X_2+\beta X_2^2\bigr)\subset\widetilde{L}.$$
We now have an explicit equation describing $L^*$. To compute $[L^*]$ we first compute its complement in $\widetilde{L}$, which is $\spec\f[X_1,X_2]/(X_1^2-\alpha X_1X_2+\beta X_2^2)\subset\widetilde{L}$. This can be split into two parts, the closed subscheme $\spec\f[X_1]/(X_1^2)$, which maps to $1$ in $\KO(\catvar_{\f})$, and its complement
$$\spec\frac{\f[X_1,X_2,1/X_2]}{\bigl(X_1^2-\alpha X_1X_2+\beta X_2^2\bigr)}
\simeq\spec\frac{\f[Y_1,Y_2,1/Y_2]}{\bigl(Y_1^2-\alpha Y_1+\beta\bigr)}.$$
Now if $\kar\f\neq2$ then $\alpha=0$ so $Y_1^2-\alpha Y_1+\beta=f(Y_1)$ and this is also true if $\kar\f=2$ for then $-\alpha=\alpha$. Hence the above expression equals $\spec\f[Y_2,1/Y_2]\times_{\f}\spec L$, and this maps to $(\bbL-1)\cdot[\spec L]\in\KO(\catvar_{\f})$. Putting this together gives $[L^*]=\bbL^2-[\spec L]\cdot\bbL+[\spec L]-1\in\KO(\catvar_{\f})$.
\end{example}
This method does not work when $L$ is a field of degree greater than $2$, the cutting and pasting then becomes to complicated. The rest of this section is devoted to giving a systematic way of computing $[L^*]$.

\subsection*{Reduction to case of lower dimension}
We now describe a method which makes it possible to reccursively compute $[L^*]$ as an element of $\artcl_{\f}[\bbL]$. This will be done in the following way. We first describe subschemes of $\widetilde{L}$, denoted  $L_1\dotsk L_n$, such that $[L^*]=\bbL^n-\sum_{i=1}^n[L_i]$. We are then reduced to compute $[L_i]$ for every $i$. To do that we find a subscheme $T_i$ of $L_i$ and a $T_i$-algebra $L'_i$ of rank $n-i$. More precisely, $T_i$ is the spectrum of a finite product of fields $\prod K_v$, where $L'_i$ equals $L_v$ on $K_v$, We then show that $L_i\simeq(L'_i/T_i)^*$ as $\f$-schemes. In Lemma \ref{90} we will show that  $(L'_i/T_i)^*\simeq\disjunion(L_v/K_v)^*$ and we are then in the situation we started with, only that the algebras have dimension less than $n$, for having computed $[(L_v/K_v)^*]\in\KO(\catvar_{K_v})$ we can find $[(L_v/K_v)^*]\in\KO(\catvar_{\f})$ with the help of the projection formula.

To do this we will need the following lemma.
\begin{lemma}\label{90}
Let $K=\prod_{v\in I}K_v$ where the $K_v$ are fields and $I$ is finite. Let $L$ be a free $K$-algebra of rank $n$, i.e., $L=\prod_{v\in I}L_v$ where, $L_v$ is a $K_v$-algebra of dimension $n$. Then $\widetilde{L/K}\simeq\disjunion_v\widetilde{L_v/K_v}$ and $(L/K)^*\simeq\disjunion_v(L_v/K_v)^*$ as $K$-schemes.
\end{lemma}
\begin{proof}
The first part follows since $\Sym(L^\vee)\simeq\prod_{v\in I}\Sym(L_v^\vee)$ as $K$-algebras. To prove that $(L/K)^*\simeq\disjunion(L_v/K_v)^*$ as $K$-schemes we prove that their functors of points are equal. Let $R$ be a $K$-algebra, i.e., $R=\prod_vR_v$ where $R_v$ is a $K_v$-algebra, possibly equal to zero. An $R$-point on $\disjunion(L_v/K_v)^*$ is a morphism $f\colon\disjunion_v\spec R_v\to\disjunion(L_v/K_v)^*$ that commutes with the structural morphisms to $\disjunion_v\spec K_v$. Since the image of $\spec R_v$ under the structural morphism is contained in $\spec K_v$ we must have $f(\spec R_v)\subset L^*_v$. Therefore $f$ is determined by a set of morphisms $\{f_v\colon\spec R_v\to L_v^*\}_{v\in I}$ where $f_i$ is a morphism of $K_v$-schemes. Hence we can identify $f$ with an element in $\prod L_v^*(R_v)$. The same is true for an $R$-point on $L^*$ for $$L^*(R)=\biggl(\Bigl(\prod L_v\Bigr)\otimes_{\prod K_v}\Bigl(\prod R_v\Bigr)\biggr)^{\times}\simeq\prod(L_v\otimes_{K_v}R_v)^\times=\prod L_v^*(R_v).$$
So by Yoneda's lemma, $L^*\simeq\disjunion L_v^*$. (This method could also have been used to prove the first part of the lemma, but there we knew the algebra representing $\widetilde{L}$ and that gave a shorter proof.)
\end{proof}
We will now give the definitions of $L_i$, $T_i$ and $L'_i$. There are two ways of doing this. The first is to construct them explicitly in much the same way as we constructed $L^*$ with the help of the norm map. The second is to just construct their images after scalar extension, as subschemes of $(\widetilde{L})_{\sepclf}$, and then use Galois descent as described in Section \ref{2}. The first method requires more work but has the advantage that it works also when the base is not a field. It is carried out in \cite{rokaeus07}. In this paper we will use the second method:

Let $L$ be the separable $\f$-algebra and let $S$ be the corresponding $\G$-set i.e., the corresponding $\fG$-algebra is $\sepcl{\f}^S$. Since, for every $\sepclf$-algebra $R$, $$\widetilde{L}_{\sepclf}(R)=(L\otimes_{\f}\sepclf)\otimes_{\sepclf}R=\sepclf^S\otimes_{\sepclf}R=(\sepclf\otimes_{\sepclf}R)^S=R^S=\bbA_{\sepclf}^S(R),$$
it follows that the $\f$-scheme $\widetilde{L}$ corresponds to the $\fG$-scheme $\bbA_{\sepclf}^S$, i.e., $\spec\sepclf[X]_{s\in S}$ where $\G$ acts on the scalars and by permuting the indeterminates. In the same way we see that the $\f$-scheme $L^*$ corresponds to the $\fG$-scheme $\multgroup^S$. (This also gives an alternative construction of $L^*$ when the base is a field, it is the $\f$-scheme corresponding to the $\fG$-scheme $\multgroup^S\subset\bbA_{\sepclf}^S$.)

We next define $L_i\subset\widetilde{L}$, where $i\in\{0,\dots,n\}$. For this, let $\cc{P}_i(S)$ be the $\G$-set of subsets of $S$ of cardinality $i$. We then define $L_i$ to be the $\f$-scheme corresponding to the $\fG$-scheme
$$\bigcup_{T\in\cc{P}_i(S)}\multgroup^{S\setminus T}\subset\bbA_{\sepclf}^S,$$
where $\multgroup^{S\setminus T}(R)$ is the group of $n$-tuples $(r_s)_{s\in S}\in\bbA^S(R)$ such that $r_s=0$ if $s\in T$ and $r_s\in R^\times$ if $r_s\notin T$. Since $$\multgroup^{S\setminus T}=\V\bigl(\{X_s\}_{s\in T}\bigr)\setminus\V\bigl(\{X_s\}_{s\notin T}\bigr)\subset\spec\sepclf[X_s]_{s\in S}=\bbA_{\sepclf}^S$$
we see that $\bigcup_{T\in\cc{P}_i(S)}\multgroup^{S\setminus T}$ is locally closed and that their union over all $i$ cover $\bbA_{\sepclf}^S$. It hence follows that the $L_i$ are locally closed and that they cover $\widetilde{L}$. Noting that $L_0=L^*$, we see that $[L^*]=\bbL^n-\sum_{i=1}^n[L_i]\in\KO(\catvar_{\f})$.

Now consider $\V(\{X_t-1\}_{t\notin T})\subset\multgroup^{S\setminus T}$. This subscheme is isomorphic to $\spec\sepclf$. Taking the union over every $T\in\cc{P}_i(S)$ we get an $\fG$-subscheme 
\begin{equation}\label{5}
\spec\sepclf^{\cc{P}_i(S)}\subset\bigcup_{T\in\cc{P}_i(S)}\multgroup^{S\setminus T}.
\end{equation}
We denote the corresponding $\f$-scheme with $T_i$. Moreover, the fact that $(T_i)_{\sepclf}=\spec\sepclf^{\cc{P}_i(S)}$ shows (in fact is equivalent to) that $T_i$ is the spectrum of a separable algebra. Hence $T_i$ corresponds to the $\G$-set $\cc{P}_i(S)$ and it is a product of separable field extensions of $\f$.

Next use the algebras $\sepclf\to\sepclf^{S\setminus T}$ for $T\in\cc{P}_i(S)$ to define a $\sepclf^{\cc{P}_i(S)}$-algebra $$\prod_{T\in\ccP_i(S)}\sepclf^{S\setminus T}$$
of rank $\absb{S\setminus T}=n-i$ in the category of $\fG$-algebras. Denote the corresponding $T_i$-algebra with $L'_i$, which is then also of rank $n-i$. We note that the map $T_i\to L'_i$ corresponds the projection $\set{(s,T)\in S\times\ccP_i(S)}{s\notin T}\to\ccP_i(S)$ in the category of $\G$-sets.

Now by Proposition \ref{90},
$$\Biggl(\prod_{T\in\ccP_i(S)}\sepclf^{S\setminus T}\bigg/\sepclf^{\cc{P}_i(S)}\Biggr)^*=\bigcup_{T\in\ccP_i(S)}\bigl(\sepclf^{S\setminus T}/\sepclf\bigr)^*$$
and since this is isomorphic to $\bigcup_{T\in\ccP_i(S)}\multgroup^{S\setminus T}$ it follows that the corresponding $\f$-schemes are isomorphic, i.e., $(L_i'/T_i)^*\simeq L_i$.

Since $T_i$ is a finite product of separable extensions of $\f$, $\prod K_v$, we must have $L'_i=\prod L_v$ where $L_v$ is a $K_v$-algebra of dimension $n-i$ (since the rank of $L'_i/T_i$ is $n-i$). Hence we see that $(L'_i/T_i)^*=\bigcup(L_v/K_v)^*$. By induction then, we may compute $[(L_v/K_v)^*]\in\KO(\catvar_{K_v})$ for each $v$ and then us the projection formula \eqref{13} to compute $[(L_i/T_i)^*]$ and hence $[L^*]\in\KO(\catvar_{\f})$.

We illustrate with an example.
\begin{example}\label{6}
Let $\f=\bbFq$ and $L=\bbF_{q^3}$. We then have
\begin{equation}\label{eq:63}
[L^*]=\bbL^3-[L_1]-[L_2]-1\in\KO(\catvar_{\f}).
\end{equation}
Let $\G:=\Gal(\overline{\f}/\f)$ and let $\fro$ be the topological generator of $\G$, the Frobenius automorphism $\alpha\mapsto\alpha^q$. Then $L$ corresponds to the $\G$-set $S:=\Hom_{\f}(L,\overline{\f})=\{1,\fro,\fro^2\}$, where we have identified $\fro$ with its restriction to $L$.

We have $\mathcal{P}_1(S)=\{\{1\},\{\fro\},\{\fro^2\}\}\simeq S$. Therefore $T_1\simeq\spec L$. Moreover, $L'_1$ corresponds to
$$\bigl\{(1,\{\fro\}),(1,\{\fro^2\}),(\fro,\{1\}),(\fro,\{\fro^2\}),(\fro^2,\{1\}),(\fro^2,\{\fro\})\bigr\}$$ and this is the union of two sets on which $\G$ acts transitive, hence it is isomorphic to $S\disjunion S$ as a $\G$-set. So $L'_1\simeq L^2$. Therefore $[(L'_1/T_1)^*]=(\bbL-1)^2\in\KO(\catvar_L)$ and hence by \eqref{13}
$$[L_1]=\ress{\f}{L}\bigl((\bbL-1)^2\bigr)=[\spec L]\cdot(\bbL-1)^2\in\KO(\catvar_{\f})$$
In the same way we find that $[L_2]=[\spec L]\cdot(\bbL-1)\in\KO(\catvar_{\f})$
Putting this into \eqref{eq:63} gives that
$$[L^*]=\bbL^3-[\spec L]\cdot\bbL^2+[\spec L]\cdot\bbL-1\in\KO(\catvar_{\f}).$$
\end{example}

\subsection*{An explicit formula}
We have now showed how to compute $[L^*]$ in principle. Evolving what we already have proved will give us an explicit formula.

To get more compact formulas we use the following notation. In the last section we worked with a fixed algebra $L/\f$ and defined $L_i$, $T_i$ and $L'_i$ with respect to this algebra. To translate the recursion into a closed formula we need to repeat this constructions. Hence we fix once and for all our base field, $\f$. Let $K$ be a finite extension field of $\f$ and let $L$ be a separable $K$-algebra. We then use $L_i(L/K)$, $T_i(L/K)$ and $L'_i(L/K)$ to denote $L_i$, $T_i$ and $L'_i$ constructed with respect to the algebra $L/K$. We have that $L_i(L/K)$ and $T_i(L/K)$ are $K$-schemes which we also can view as $\f$-schemes by restriction of scalars. Similarly, $L'_i(L/K)$ is a $T_i(L/K)$-algebra which we also may view as a $\f$-algebra.

We next wants to generalize these definitions to the case when $K$ is a finite separable $\f$-algebra. Recall that this is already done for $(L/K)^*$. However, since the definitions of $L_i$, $T_i$ and $L'_i$ use Galois decent their constructions cannot be generalized directly. We instead do the following.
\begin{definition}\label{96}
Let $K$ be a finite separable $\f$-algebra and $L$ a finite separable $K$-algebra, so $K=\prod_vK_v$ where $K_v$ are separable extension fields of $\f$ and $L=\prod_vL_v$ where $L_v$ is a separable $K_v$-algebra. Define
$$L_i(L/K):=\disjunion_vL_i(L_v/K_v).$$
Furthermore, define
$$T_i(L/K):=\disjunion_vT_i(L_v/K_v)$$
and define $L'_i(L/K)$ to be the $T_i(L/K)$-algebra which is $L'_i(L_v/K_v)$ on $T_i(L_v/K_v)$.
\end{definition}
As in this definition, let $K$ be a finite separable $\f$-algebra and $L$ a finite separable $K$-algebra of rank $n$, so $K=\prod_vK_v$ where $K_v$ are separable extension fields of $\f$ and $L=\prod_vL_v$ where $L_v$ is a separable $K_v$-algebra of dimension $n$ as a vector space over $K_v$. Using Definition \ref{96}, Proposition \ref{90} and the stratification of $\widetilde{L}$ over a field given in the preceding section gives the following:
\begin{lemma}
We have that
$$\widetilde{L/K}=(L/K)^*\cup\bigcup_{i=1}^{n-1}L_i(L/K)\cup\spec K$$
where the union is disjoint and open. Hence,
$$[(L/K)^*]=[\spec K]\cdot\bbL^n-\sum_{i=1}^{n-1}[L_i(L/K)]-[\spec K]\in\KO(\catvar_{\f}).$$
\end{lemma}
From the preceding section we know that $L'(L_v/K_v)/T_i(L_v/K_v)$ has rank $n-i$ and that $L_i(L_v/K_v)\simeq\bigl(L'_i(L_v/K_v)/T_i(L_v/K_v)\bigr)^*$. Taking the union over every $v$ and using Proposition \ref{90} we get the following.
\begin{lemma}\label{166}
The algebra $L'_i(L/K)/T_i(L/K)$ has rank $n-i$, and
$$L_i(L/K)\simeq\bigl(L'_i(L/K)/T_i(L/K)\bigr)^*$$
as $\f$-schemes.
\end{lemma}

For the rest of this section, we fix a field $\f$ and a separable $\f$-algebra $L$ of dimension $n$. Notation: Given a sequence of positive integers $i_1\dotsk i_q$. Construct the algebra $L'_{i_1}/T_{i_1}=L'_{i_1}(L/\f)/T_{i_1}(L/\f)$. Define the algebra $L'_{i_2,i_1}/T_{i_2,i_1}$ as $L'_{i_2}(L'_{i_1}/T_{i_1})/T_{i_2}(L'_{i_1}/T_{i_1})$ and define inductively $L'_{i_{r+1}\dotsk i_1}/T_{i_{r+1}\dotsk i_1}$ as
$$L'_{i_{r+1}}(L'_{i_r\dotsk i_1}/T_{i_r\dotsk i_1})/T_{i_{r+1}}(L'_{i_r\dotsk i_1}/T_{i_r\dotsk i_1}).$$
Inductively we get that the rank of $L'_{i_{r}\dotsk i_1}/T_{i_{r}\dotsk i_1}$ is $n-(i_1+\dots+i_r)$. Hence, as a corollary to the preceding lemmas we get the following.
\begin{lemma}\label{104}
Let $\alpha=(i_r\dotsk i_1)$ where $\sum_{s=1}^ri_s=i$. Then 
\begin{equation*}
\bigl[(L'_\alpha/T_\alpha)^*\bigr]=[T_\alpha]\cdot\bbL^{n-i}-\sum_{j=1}^{n-i-1}\bigl[(L'_{j,\alpha}/T_{j,\alpha})^*\bigr]-[T_\alpha]\in\KO(\catvar_{\f}).
\end{equation*}
\end{lemma}

We are now ready to prove the main theorem of this subsection.
\begin{theorem}\label{109}
With the same notation as above we have
$$[L^*]=\bbL^n+a_1\bbL^{n-1}+\dots+a_{n-1}\bbL+a_n$$
where 
\begin{equation*}
a_j=\sum_{r=1}^{j}(-1)^r\sum_{\substack{(i_1\dotsk i_{r}):\\i_1+\dots+i_{r}=j\\i_s\geq1}}[T_{i_{r}\dotsk i_1}]
\end{equation*}
for $j=1\dotsk n$.
\end{theorem}
\begin{proof}
We evaluate $[L^*]$ in $n$ steps, using Lemma \ref{104}. In the first step we write
$$[(L/\f)^*]=\bbL^n-[(L'_1/T_1)^*]-\dots-[(L'_{n-1}/T_{n-1})^*]-1$$
so we get the contribution $\bbL^n-1$. We then evaluate the remaining terms, using lemma \ref{104}, so in step two we get a sum consisting of two parts. First, $\bigl[(L'_{i_2,i_1}/T_{i_2,i_1})^*\bigr]$ shows up with sign $(-1)^2$, for $2\leq i_2+i_1<n$ (we always have $i_s\geq1$). This is the terms that we will take care of in step three. The second part of the sum contributes to our formula. It consists of the terms 
\begin{align*}
(-1)^2\bigl(-[T_j]\cdot\bbL^{n-j}+[T_j]\bigr) && 1\leq j<n.
\end{align*}
Continuing in this way we find that in step $r$ we get a sum consisting of two parts. Firstly, every term of the form $\bigl[(L'_{i_r\dotsk i_1}/T_{i_r\dotsk i_1})^*\bigr]$ with coefficient $(-1)^r$, for $\sum_{s=1}^ri_s<n$. This part is taken care of in step $r+1$. And secondly we get a contribution to our formula consisting of
\begin{align*}
(-1)^r\bigl(-[T_{i_{r-1}\dotsk i_1}]\cdot\bbL^{n-j}+[T_{i_{r-1}\dotsk i_1}]\bigr) && r-1\leq j<n
\end{align*}
for every $r-1$-tuple $(i_{r-1}\dotsk i_1)$ such that $\sum_{s=1}^{r-1}i_s=j$. This process ends in step $n$.

Collecting terms we now see that if $1\leq j\leq n-1$ then the coefficient in front of $\bbL^{n-j}$ becomes
$$\sum_{r=2}^{j+1}(-1)^{r+1}\sum_{\substack{(i_1\dotsk i_{r-1}):\\i_1+\dots+i_{r-1}=j\\i_s\geq1}}[T_{i_{r-1}\dotsk i_1}].$$
This equals 
\begin{equation}\label{105}
\sum_{r=1}^{j}(-1)^r\sum_{\substack{(i_1\dotsk i_{r}):\\i_1+\dots+i_{r}=j\\i_s\geq1}}[T_{i_{r}\dotsk i_1}].
\end{equation}

A separate computation, using that $[T_n]=1 $ and that if $1\leq\sum_{s=1}^{r-1}i_s=j<n$ then $T_{n-j,i_{r-1}\dotsk i_1}=T_{i_{r-1}\dotsk i_1}$, shows that formula \eqref{105} holds also for the constant coefficient, when $j=n$.
\end{proof}

\subsection*{The formula expressed using the Burnside ring}

The formula in the preceding secton is not suitable for computations. In this section we make it, by computing the object in the Burnside ring that maps to $[T_\alpha]$. We begin with some notations.
\begin{definition}
Let $G$ be a profinite group. Given a $G$-set $S$ of cardinality $n$ and a positive integer $r$. Moreover, let $(i_1\dotsk i_r)$ be an $r$-tuple of positive integers such that $i_1+\dots+i_r\leq n$. Then $\mathcal{P}_{i_r\dotsk i_1}(S)$ is the $G$-set of $r$-tuples $(S_r\dotsk S_1)$ where $S_{j}$ is a subset of $S$ of cardinality $i_j$ and the $S_j$ are pairwise disjoint. In particular $\mathcal{P}_i(S)$ has the same meaning as before (up to isomorphism).
\end{definition}
\begin{lemma}\label{lemma:64}
Let $\f$ be a field and $K$ a separable $\f$-algebra of dimension $t$. Let $L$ be a separable $K$-algebra of rank $n$. Let $\G:=\Gal(\f^s/\f)$ and let $K$ and $L$ correspond to $T$ respectively $S$ as $\G$-sets. Write $T=\Hom_{\f}(K,\f^s)=\{\tau_1\dotsk\tau_t\}$. Let $S_j$ be the inverse image of $\tau_j$ under the map $S\rightarrow T$ corresponding to $K\rightarrow L$. Then $T_i(L/K)$ corresponds to the $\G$-set
$$\bigcup_{j=1}^t\mathcal{P}_i(S_j)$$
and $L'_i(L/K)$ corresponds to
$$\bset{(f,U)\in\bigcup_{j=1}^tS_j\times\mathcal{P}_i(S_j)}{f\notin U}$$
\end{lemma}
\begin{proof}
Suppose first that $K$ is a field. According to \eqref{5}, $T_i(L/K)$ corresponds to $\mathcal{P}_i\bigl(\Hom_K(L,\f^s)\bigr)$ as a $\Gal(\f^s/K)$-set. Hence by Proposition \ref{Pp:60} it corresponds to
$$\G\times_{\Gal(\f^s/K)}\mathcal{P}_i\bigl(\Hom_K(L,\f^s)\bigr)$$
as a $\G$-set, with the $\G$-action given in that proposition. Since we assumed that $K$ is a field we may write $T$ as $\{\tau_1\vert_K\dotsk\tau_t\vert_K\}$, where $\tau_j\in\G$, and this in turn can be identified with a system of coset representatives of $\G/\Gal(\f^s/K)$. We hence want to show that we have an isomorphism of $\G$-sets,
$$
\phi\colon T\times\mathcal{P}_i\bigl(\Hom_K(L,\f^s)\bigr)\rightarrow\bigcup_{j=1}^t\mathcal{P}_i(S_j).
$$
To construct this, define $\phi$ as $(\tau_j\vert_K,U)\mapsto\tau_jU$. (Note that $\tau_j$ have to be fixed for every $j$, if we replace it with $\tau'_j$ such that $\tau_j\vert_K=\tau'_j\vert_K$ we may get another $\phi$.)  First $\phi$ is well defined because every element in $U$ fixes $K$, so every element of $\tau_jU$ is in $S_j$, the inverse image of $\tau_j\vert_K$ in $S$. Hence $\phi(\tau_j\vert_K,U)\in\mathcal{P}_i(S_j)$. It is also $\G$-equivariant, because if $\sigma\in\G$ is such that $\sigma\tau_j=\tau_l\tau'$, where $\tau'\in\Gal(\f^s/K)$, then
$$\phi\bigl(\sigma(\tau_j\vert_K,U)\bigr)=\phi(\tau_l,\tau'U)=\tau_l\tau'U$$
and
$$\sigma\phi(\tau_j\vert_K,U)=\sigma(\tau_jU)=\sigma\tau_jU=\tau_l\tau'U.$$
Next $\phi$ is injective: If $\phi(\tau_j\vert_K,U)=\phi(\tau_l\vert_K,U')$ then they both must be in $\mathcal{P}_i(S_j)$, so $l=j$. Hence $\tau_jU=\tau_jU'$ and since $\tau_j$ is an isomorphism, $U=U'$. So $\phi$ is an injective morphism between two $\G$-sets of cardinality $t\cdot\tbinom{n}{i}$, hence an isomorphism.

For the general case when $K$ is a separable $\f$-algebra of dimension $t$, note that we can identify $T$ with
$$\Disjunion_v\Hom_{\f}(K_v,\f^s)$$
where $K=\prod_vK_v$, by sending $f\in\Hom_{\f}(K_{v_0},\f^s)$ to $(\alpha_v)\mapsto f(\alpha_{v_0})\in T$. Denote the map $S\rightarrow T$ by $\pi$. We have that $T_i(L/K)=\disjunion_vT_i(L_v/K_v)$. This corresponds to the $\G$-set
$$\Disjunion_v\bigcup_{\tau\in\Hom_{\f}(K_v,\f^s)}\mathcal{P}_i(\pi^{-1}\tau)=\bigcup_{\tau\in T}\mathcal{P}_i(\pi^{-1}\tau)=\bigcup_{j=1}^t\mathcal{P}_i(S_j)$$

As for $L'_i(L/K)$, assume first that $K$ is a field. As a $\Gal(\f^s/K)$-set, $L'_i(L/K)$ corresponds to
$$M:=\set{(f,U)\in\Hom_K(L,\f^s)\times\mathcal{P}_i\bigl(\Hom_K(L,\f^s)\bigr)}{f\notin U},$$
hence it corresponds to $T\times M$ as a $\G$-set. Define a map
$$T\times M\rightarrow\bset{(f,U)\in\bigcup_{j=1}^tS_j\times\mathcal{P}_i(S_j)}{f\notin U}$$
by
$$\bigl(\tau_j\vert_K,(f,U)\bigr)\mapsto(\tau_j\circ f,\tau_jU).$$
As above one shows that this is an isomorphism of $\G$-sets. The case when $K$ is an arbitrary separable $\f$-algebra is handled in the same way as $T_i$.
\end{proof}

\begin{proposition}\label{110}
Let $\alpha=(i_r\dotsk i_1)$ be an $r$-tuple of positive integers such that $i_1+\dots+i_r=i$ where $1\leq i\leq n$. The algebra $L'_\alpha/T_\alpha$ in the category of $\f$-algebras corresponds to the $\G$-sets
$$\bbset{\bigl(s,(S_r\dotsk S_1)\bigr)\in S\times\mathcal{P}_\alpha(S)}{s\notin\cup_{t=1}^rS_t}$$
and $\mathcal{P}_\alpha(S)$ together with the projection morphism.
\end{proposition}
\begin{proof}
By construction the proposition holds for $r=1$. Suppose the formula has been proved for $r$. We have $T_{i_{r+1},i_r\dotsk i_1}=T_{i_{r+1}}(L'_{i_r\dotsk i_1}/T_{i_r\dotsk i_1})$. By the induction hypothesis and lemma \ref{lemma:64} this corresponds to
\begin{equation*}
\bigcup_{(S_r\dotsk S_1)\in\mathcal{P}_{i_r\dotsk i_1}(S)}\mathcal{P}_{i_{r+1}}\Bigl(\bbset{\bigl(s,(S_r\dotsk S_1)\bigr)}{s\notin\cup_{t=1}^rS_t}\Bigr)
\end{equation*}
which is isomorphic to
\begin{equation*}
\bigcup_{(S_r\dotsk S_1)\in\mathcal{P}_{i_r\dotsk i_1}(S)}\bbset{\bigl(\{s_1\dotsk s_{i_{r+1}}\},S_r\dotsk S_1\bigr)}{s_{i_t}\notin\cup_{t=1}^rS_t}
\end{equation*}
and this in turn is equal to $\mathcal{P}_{i_{r+1},i_r\dotsk i_1}(S)$.

And $L'_{i_{r+1},i_r\dotsk i_1}$ corresponds to the set of pairs $(f,U)$ in
\begin{equation*}
\bigcup_{(S_r\dotsk S_1)\in\mathcal{P}_{i_r\dotsk i_1}(S)}\bbset{\bigl(s,(S_r\dotsk S_1)\bigr)}{s\notin\cup_{t=1}^rS_t}\times\mathcal{P}_{i_{r+1}}\Bigl(\bbset{\bigl(s,(S_r\dotsk S_1)\bigr)}{s\notin\cup_{t=1}^rS_t}\Bigr)
\end{equation*}
such that $f\notin U$. This set is isomorphic to
\begin{equation*}
\bigcup_{(S_r\dotsk S_1)\in\mathcal{P}_{i_r\dotsk i_1}(S)}\bbset{\bigl(s,(S_{r+1},S_r\dotsk S_1)\bigr)\in S\times\cc{P}_{i_{r+1},i_r\dotsk i_1}(S)}{s\notin\cup_{t=1}^{r+1}S_t}
\end{equation*}
which equals
\begin{equation*}
\bbset{\bigl(s,(S_{r+1},S_r\dotsk S_1)\bigr)\in S\times\mathcal{P}_{i_{r+1}\dotsk i_1}(S)}{s\notin\cup_{t=1}^{r+1}S_t}.\qedhere
\end{equation*}
\end{proof}

We are now ready to give our first closed formula for $[L^*]$. It follows from Theorem \ref{109} and Proposition \ref{110}.
\begin{theorem}\label{82}
Let $L$ be a $\f$-algebra of dimension $n$ and $S$ a $\G$-set such that $\art_{\f}\bigl([S]\bigr)=[\spec L]$. Define
\begin{equation*}
\rho_i(S)=\sum_{t=1}^{i}\sum_{\substack{(i_1\dotsk i_t):\\i_1+\dots+i_t=i\\i_s\geq 1}}(-1)^t[\mathcal{P}_{i_t\dotsk i_1}(S)]\in\bur(\G).
\end{equation*}
Then
$$[L^*]=\bbL^n+a_1\cdot\bbL^{n-1}+\dots+a_{n-1}\cdot\bbL+a_n\in\KO(\catvar_{\f})$$
where $a_i=\art_{\f}(\rho_i(S))$.
\end{theorem}

\subsection*{The universal nature of the formula}

Fix a field $\f$ with absolute Galois group $\G$. Also, fix a separable $\f$-algebra $L$ of dimension $n$ corresponding to the $\G$-set $S$. Define a homomorphism $\phi\colon\G\to\per_n$ as the composition of $\G\to\aut(S)$ with an isomorphism $\aut(S)\to\per_n$. Let $\res_{\G}^{\per_n}\colon\bur(\per_n)\to\bur(\G)$ be the restriction map with respect to $\phi$. Then $\res_{\G}^{\per_n}$ is independent of the chosen isomorphism $\aut(S)\to\per_n$. We have that
$$\res_{\G}^{\per_n}\bigl(\{1\dotsk n\}\bigr)=[S]\in\bur(\G).$$
Using the notation that $\ccP^{(n)}_\alpha:=\ccP_\alpha(\{1\dotsk n\})$ we also have $\res_{\G}^{\per_n}(\cc{P}_\alpha)=[\cc{P}_\alpha(S)]$. Define
\begin{equation}\label{22}
\rho_i^{(n)}=\sum_{t=1}^{i}\sum_{\substack{(i_1\dotsk i_t):\\i_1+\dots+i_t=i\\i_s\geq 1}}(-1)^t[\mathcal{P}^{(n)}_{i_t\dotsk i_1}]\in\bur(\per_n).
\end{equation}
Then $\res_{\G}^{\per_n}(\rho_i^{(n)})=\rho_i(S)$. This discussion gives the following formulation of Theorem \ref{82}.
\begin{theorem}\label{111}
Fix a positive integer $n$. The $\rho_i^{(n)}\in\bur(\per_n)$, defined in \eqref{22}, are universal in the sense that for every field $\f$ with absolute Galois group $\G$ and every separable $\f$-algebra of dimension $n$,
$$[L^*]=\bbL^n+a_1\cdot\bbL^{n-1}+\dots+a_{n-1}\cdot\bbL+a_n\in\KO(\catvar_{\f}),$$
where $a_i=\art_{\f}\circ\res_{\G}^{\per_n}(\rho_i^{(n)})$.
\end{theorem}
We illustrate with an example.
\begin{example}
Let $L/\f=\bbF_{q^4}/\bbFq$. Since $\G$ is generated by the Frobenius map $\fro$ we can identify $S$, the $\G$-set corresponding to $L$, with $\{1,\fro,\fro^2,\fro^3\}$. We have
\begin{equation*}
\cc{P}_2(S)=\bigl\{\{1,\fro\},\{\fro,\fro^2\},\{\fro^2,\fro^3\},\{1,\fro^3\}\bigr\}\disjunion\bigl\{\{1,\fro^2\},\{\fro,\fro^3\}\bigr\}.
\end{equation*}
The first of these sets is isomorphic to $S$. The second is transitive of cardinality $2$ so it corresponds to a field extension of $\f$ of degree $2$, i.e., $\bbF_{q^2}$. Reasoning in this way and using theorem \ref{111} we find that
\begin{equation*}
[L^*]=\bbL^4-[\spec\bbF_{q^4}]\cdot\bbL^3+\bigl(2[\spec\bbF_{q^4}]-[\spec\bbF_{q^2}]\bigr)\cdot\bbL^2-[\spec\bbF_{q^4}]\cdot\bbL+[\spec\bbF_{q^2}]-1.
\end{equation*}

If instead $L/\f=\bbF_{q^2}\times\bbF_{q^2}/\bbFq$ then $S=\{e_1,\fro e_1\}\disjunion\{e_2,\fro e_2\}$ where $e_1$ and $e_2$ are the projection maps. We then get, for example,
\begin{equation*}
\cc{P}_2(S)=\bigl\{\{e_1,\fro e_1\}\bigr\}\disjunion\bigl\{\{e_2,\fro e_2\}\bigr\}\disjunion\bigl\{\{e_1,e_2\},\{\fro e_1,\fro e_2\}\bigr\}\disjunion\bigl\{\{e_1,\fro e_2\},\{\fro e_1,e_2\}\bigr\}.
\end{equation*}
This kind of computations gives that
\begin{equation}\label{23}
[L^*]=\bbL^4-2[\spec\bbF_{q^2}]\cdot\bbL^3+\bigl(4[\spec\bbF_{q^2}]-2\bigr)\cdot\bbL^2-2[\spec\bbF_{q^2}]\cdot\bbL+1.
\end{equation}
Note  however that since $L^*(R)=(\bbF_{q^2}\otimes_{\f}R)^\times\times(\bbF_{q^2}\otimes_{\f} R)^\times=(\bbF_{q^2}^*\times_{\f}\bbF_{q^2}^*)(R)$, Yoneda's lemma shows that $L^*\simeq\bbF_{q^2}^*\times_{\f}\bbF_{q^2}^*$, so \eqref{23} could also have been obtained by squaring the expression for $[\bbF_{q^2}^*]$ given in Example \ref{9}.
\end{example}

\section{The class of the torus in terms of the $\lambda$-operations}\label{12}

In \cite{rokaeusBur} we study the $\lambda$-structure on $\bur(\per_n)$. The main result is that there is a formula for $\lambda^i$ namely that, given $n$, for $i=1\dotsk n$ we have
\begin{equation}\label{14}
\lambda^i\bigl(\{1\dotsk n\}\bigr)=(-1)^i\sum_{t=1}^{i}\sum_{\substack{(i_1\dotsk i_t):\\i_1+\dots+i_t=i\\i_s\geq 1}}(-1)^t\bigl[\mathcal{P}_{i_t\dotsk i_1}^{(n)}\bigr]\in\bur(\per_n).
\end{equation}
From this and Theorem \ref{82} the following theorem is immediate.
\begin{theorem}\label{119}
Let $\rho_i^{(n)}$ be the elements defined in Theorem \ref{111}, i.e., the elements in $\bur(\per_n)$ describing $[L^*]\in\KO(\catvar_{\f})$ for every separable, $n$-dimensional algebra $\f\to L$. Then
$$\rho_i^{(n)}=(-1)^i\cdot\lambda^i\bigl(\{1\dotsk n\}\bigr)\in\bur(\per_n).$$
\end{theorem}
As a corollary, \eqref{3} follows:
\begin{proof}[Proof of \eqref{3}]
Let $L$ correspond to the $\G$-set $S$. From Theorem \ref{111} we know that the coefficient in front of $\bbL^{n-i}$ is $\art_{\f}\circ\res_{\G}^{\per_n}(\rho_i^{(n)})$ and by theorem \ref{119} this equals $\art_{\f}\bigl((-1)^i\cdot\lambda^i(S))$. Since $\art_{\f}$ is a $\lambda$-morphism that maps $[S]$ to $[\spec L]$ the result follows.
\end{proof}
We now give an alternative proof of the equality $\rho_i^{(n)}=(-1)^i\lambda^i\bigl(\{1\dotsk n\}\bigr)\in\bur(\per_n)$. This proof is based on  point counting over finite fields. We consider it motivated to include since one purpose of this article is to give examples of different computation techniques in $\KO(\catvar_{\f})$.

This proof do not use the universal formula \eqref{14} from \cite{rokaeusBur}. We do however need the following results proved in \cite{rokaeusBur}: Write $\ccP_{\mu}$ for the $\per_n$-set $\ccP_\mu\bigl(\{1\dotsk n\}\bigr)$. We define the Schur subring $\schur_n\subset\bur(\per_n)$ to be the subgroup generated by $\{[\ccP_\mu]\}_{\mu\vdash n}$. (Here $\mu\vdash n$ means that $\mu$ is a partition of $n$.) This is closed under multiplication, hence really a ring. It is not a $\lambda$-ring since it is not closed under the $\lambda$-operations. However, $\lambda^i(\{1\dotsk n\})\in\schur_n$ for every $i$. Moreover the restriction of $\btr\colon\bur(\per_n)\to\rep_{\bbQ}(\per_n)$ to $\schur_n$ is injective.

Also, we do not need to know the explicit description of the universal elements $\rho_i^{(n)}$ given in Theorem \ref{111}. We do however need their existence and that they lie in $\schur_n$, and once that is proved it is not such a long step to describe the elements. If we instead use Theorem \ref{111} the below proof of Theorem \ref{119} gives a (very \emph{ad hoc}) proof of \eqref{14}.

The setting for the proof is as follows. Let $L$ be a separable $\bbFq$-algebra of dimension $n$, corresponding to the $\G:=\Gal\bigl(\algcl{\bb{F}}_q/\bbFq\bigr)$-set $S$. Choose an isomorphism $\aut(S)\to\per_n$ and compose it with the homomorphism $\G\to\aut(S)$ to get a homomorphism $\phi\colon\G\to\per_n$. Let $\fro$ be the topological generator of $\G$ and define $\sigma:=\phi(\fro)\in\per_n$. Let $\res_{\G}^{\per_n}$ denote the restriction maps with respect to $\phi$ for Burnside as well as representation rings.

Recall that we use $\ifc{\sigma}$ and $\ifc{\fro}$ for the character homomorphisms with respect to $\sigma$ and $\fro$ respectively. One sees that $\ifc{\sigma}=\ifc{\fro}\circ\res_{\G}^{\per_n}$. The corresponding map on the Grothendieck ring of varieties is the point counting homomorphism $\con{q}$ defined by $X\mapsto\abs{X(\bb{F}_{q})}$, and by \eqref{20} we have $\ifc{\fro}\circ\btr=\con{q}\circ\art_{\bbFq}$. So with this notation the following diagram is commutative:
\begin{equation}\label{11}
\xymatrix{
\bur(\per_n)  \ar[r]^{\res_{\G}^{\per_n}}\ar[d]^{\btr}          & \bur(\G)  \ar[d]^{\btr}\ar[rd]^{\art_{\bbFq}}   \\
\rep_{\bbQ}(\per_ n) \ar[r]^{\res_{\G}^{\per_n}}\ar[dr]_{\ifc{\sigma}}  & \rep_{\bbQ}(\G) \ar[d]^{\ifc{\fro}}  &   \KO(\catvar_{\bbFq})  \ar[ld]^{\con{q}}    \\
                                    & \bbZ
}
\end{equation}
Using this we are now ready to give the alternative proof.
\begin{proof}[Proof of Theorem \ref{119}]
Fix a positive integer $n$. We write $\ell_i:=\lambda^i(\{1\dotsk n\})$ and we want to prove that $\rho_i=(-1)^i\ell_i\in\bur(\per_n)$. Since they both lie in $\schur_n$ it suffices to show that if $R$ is a set of representatives of the conjugacy classes of $\per_n$ then for every $\sigma\in R$,
$$\ifc{\sigma}\btr(\rho_i)=(-1)^i\ifc{\sigma}\btr(\ell_i)\in\bbZ.$$
(Since $\btr$ is injective on $\schur_n$ and $\prod_{\sigma\in R}\ifc{\sigma}$ is injective.) We do this simultaneously for $i=0\dotsk n$ by showing that 
\begin{equation}\label{123}
\sum_{i=0}^n\ifc{\sigma}\btr(\rho_i)X^{n-i}=\sum_{i=0}^n(-1)^i\ifc{\sigma}\btr(\ell_i)X^{n-i}\in\bbZ[X]
\end{equation}
for every $\sigma\in R$.

From now on, fix $\sigma\in R$. Let $q$ be an arbitrary prime power, let $\f=\bbFq$ and let $\G:=\Gal(\overline{\f}/\f)$. Choose a $\G$-set $S$ and an isomorphism $\phi\colon\aut(S)\to\per_n$ such that $\F\mapsto\sigma$, where $\F$ is the topological generator of $\G$, the Frobenius automorphism. (Equivalently, let $S=\disjunion_{1\leq j\leq m} T_j$ such that $T_j$ is a transitive $\G$-set of cardinality $n_j$, where $\sigma$ has cycle-type $(n_1\dotsk n_m)$. Such an $S$ always exists for by Galois theory it comes from $L=\prod_{j=1}^mK_j$ where $K_j$ is a degree $n_j$ field extension of $\f$, i.e., $K_j=\bb{F}_{q^{n_j}}$.) Construct $\res_{\G}^{\per_n}$ with respect to $\phi$. In what follows we write $\ell_i(S)$ and $\rho_i(S)$ for $\res_{\G}^{\per_n}(\ell_i)$ and $\res_{\G}^{\per_n}(\rho_i)$ respectively.

We begin by computing the right hand side of \eqref{123} in terms of $(n_1\dotsk n_m)$. Let $f$ be an endomorphism of the vector space $V$ of dimension $n$. From linear algebra we know the following expression for the characteristic polynomial of $f$:
\begin{equation*}
\det(X\cdot E_n-f)=\sum_{i=0}^n(-1)^i\tr(\ext^if)X^{n-i}.
\end{equation*}
Putting $f=\F$ gives
\begin{equation}\label{121}
\det(X\cdot E_n-\F)=\sum_{i=0}^n(-1)^i\chi_{\ext^i\bbQ[S]}(\F)X^{n-i}.
\end{equation}
Since $\btr(\ell_i(S))=\bigl[\ext^i\bbQ[S]\bigr]\in\rep_{\bbQ}(\G)$ we have that $\ifc{\F}\btr(\ell_i(S))=\chi_{\ext^i\bbQ[S]}(\F)$, hence, by the commutativity of \eqref{11}, the right hand side of \eqref{121} equals
\begin{equation*}
\sum_{i=0}^n(-1)^i\ifc{\sigma}\btr(\ell_i)X^{n-i}.
\end{equation*}
As for the left hand side of \eqref{121}, since $S$ is a union of transitive $\G$-sets $T_j$ we have $\bbQ[S]=\oplus_{j=1}^m\bbQ[T_j]$ where $\bbQ[T_j]$ is irreducible, hence the matrix for $F$ is of the form 
\begin{equation*}
\begin{pmatrix}
M_1          &     &        &\text{\Large{0}}      \\
             & M_2 &        &                 \\
             &     &\ddots  &                 \\
\text{\Large{0}}    &     &        &  M_m
\end{pmatrix}
\end{equation*}
where $M_j$ is a transitive $n_j\times n_j$ permutation matrix. Since the characteristic polynomial of such a matrix is $X^{n_j}-1$ it follows that $\det(XE_n-F)=\prod_{j=1}^m\det(XE_{n_j}-M_j)=\prod_{j=1}^m(X^{n_j}-1)$. From \eqref{121} we therefore get
\begin{equation}\label{126}
\prod_{j=1}^m(X^{n_j}-1)=\sum_{i=0}^n(-1)^i\ifc{\sigma}\btr(\ell_i)X^{n-i}.
\end{equation}

We next compute the left hand side of \eqref{123}. By the definition of the $\rho_i$ we have 
$$[L^*]=\sum_{i=0}^n\art\bigl(\rho_i(S)\bigr)\bbL^{n-i}\in\KO(\catvar_{\f}).$$
Applying $\con{q}$ to this gives
\begin{equation}\label{124}
\abs{L^*(\f)}=\sum_{i=0}^n\con{q}\art(\rho_i(S))\cdot q^{n-i}.
\end{equation}
By the commutativity of \eqref{11}, $\con{q}\art(\rho_i(S))=\ifc{\sigma}\btr(\rho_i)$, so the right hand side of \eqref{124} equals
$$\sum_{i=0}^n\ifc{\sigma}\btr(\rho_i)q^{n-i}.$$
On the other hand, since we saw that $L=\prod_{j=1}^m\bb{F}_{q^{n_j}}$ we have $L^*(\f)=L^\times=\prod_{j=1}^m\bb{F}_{q^{n_j}}^\times$ so $\abs{L^*(\f)}=\prod_{j=1}^m(q^{n_j}-1)$.
Hence \eqref{124} says that
$$\prod_{j=1}^m(q^{n_j}-1)=\sum_{i=0}^n\ifc{\sigma}\btr(\rho_i)q^{n-i}.$$
Since $q$ is an arbitrary prime power it follows that
\begin{equation}\label{125}
\prod_{j=1}^m(X^{n_j}-1)=\sum_{i=0}^n\ifc{\sigma}\btr(\rho_i)X^{n-i}.
\end{equation}

Comparing \eqref{126} to \eqref{125} now gives \eqref{123}.
\end{proof}

\bibliography{../artlic/bibartlic}
\bibliographystyle{amsalpha}

\end{document}